\documentclass[11 pt]{amsart}
\usepackage{amssymb}
\usepackage{amsmath}
\usepackage{amsfonts}
\usepackage{graphicx}
\usepackage{amsthm}
\usepackage{enumerate}
\usepackage[mathscr]{eucal}
\usepackage{verbatim}

 \theoremstyle{plain}
\newtheorem{theorem}{Theorem}[section]
\newtheorem{lemma}[theorem]{Lemma}

\numberwithin{equation}{section}
\theoremstyle{plain}

\theoremstyle{remark}

\def\bbR{{\mathbb {R}}}

\begin{document}

\date{November, 2009}

\title
{Convolution with measures on flat curves in low dimensions}


\author
{Daniel M. Oberlin}

\address
{D. M.  Oberlin \\
Department of Mathematics \\ Florida State University \\
 Tallahassee, FL 32306}
\email{oberlin@math.fsu.edu}

\begin{abstract}
We prove $L^p \rightarrow L^q$ convolution estimates for the affine arclength measure on certain flat curves  
in $\bbR ^d$ when $d\in \{2,3,4\}$.
\end{abstract}

\maketitle

\section{Introduction}\label{intro}

Let $\gamma$ be a curve in $\bbR^d$ given by 

\begin{equation}\label{ineq0}
\gamma (t) =\Big(t,\frac{t^2}{2} ,\dots ,\frac{t^{d-1}}{(d-1)!},\phi (t)\Big)
\end{equation}

\noindent where $\phi\in C^{(d)}(a,b)$, where $\phi ^{(j)}(t)> 0$ for $t\in (a,b)$ and
$j=0,1,2,\dots ,d$, and where $\phi^{(d)}$ is nondecreasing. Such curves are termed {\it simple} in \cite{DM1}.
We are interested in the possibility of proving $L^p \rightarrow L^q$ convolution estimates 
for the affine arclength measure $\lambda$ on \eqref{ineq0}, 
given by $d\lambda =\phi^{(d)} (t)^{2/(d^2+d)}dt$.  We begin by recalling a theorem from \cite{O1}.  
(In this note, $|E|$ will stand for the Lebesgue measure of $E$.)

\begin{theorem} \label{theorem1} Suppose $d=2$. The inequality

\begin{equation*}
\|\lambda\ast\chi_E \|_3 \leq (12)^{1/3} \, |E|^{2/3}
\end{equation*}

\noindent holds for all measurable $E\subset\bbR^2$.
\end{theorem}
 
 \noindent Theorem \ref{theorem1} is equivalent to a weak-type $(3/2,3)$ estimate for the operator given by 
 convolution with $\lambda$, an estimate which is uniform over the class of measures $\lambda$ described above.  Here are two  questions which are raised by Theorem \ref{theorem1}:
 
 \noindent (i) is there an analogous strong-type estimate, and
 
 \noindent (ii) are there analogs of Theorem \ref{theorem1} if $d>2$? 

\noindent Having no idea how to attack these interesting questions in their natural generality, we follow 
the usual practice of asking 
what can be said along such lines by imposing additional hypotheses on $\phi$. The requirement 
\begin{equation}\label{ineq1}
\Big(\prod_{j=1}^n \phi^{(d)}(s_j)\Big)^{1/n} \le A \,  \phi^{(d)}
\big(\tfrac{s_1+\dots+s_n}{n}\big),
\end{equation}
to hold for $s_j \in (a,b)$, was used with $n=d$ in \cite{BOS2} to obtain 
Fourier restriction estimates for curves \eqref{ineq0}.
It is obvious that if $\beta\ge d$ then condition \eqref{ineq1}
holds with $A=1$  for $\phi(t)=t^\beta$ on the interval
$(0,\infty)$. Moreover, as was observed in \cite{BOS2},  if we define
$\phi_0(t)=t^{\beta}$ for some $\beta >d$  and then define 
$$
\phi_j(t)=
\int_0^t {(t-u)^{d-1}}
\exp\Big(-\tfrac{1}{\phi_{j-1}^{(d)}(u)}\Big) du
$$
for $j\ge1$, each of the functions 
$\phi_j$ satisfies \eqref{ineq1} with $A=1$ on
$(0,\infty)$. This
yields  a  sequence of functions which are
progressively flatter at the origin. 
(See \S 4 of \cite{BOS2} for other examples of flat functions satisfying \eqref{ineq1}.)   
In this note we will assume the $n=2$ version of \eqref{ineq1} which,
with $\omega\doteq (\phi^{(d)})^{2/(d^2+d)}$, we write as
\begin{equation}\label{ineq2}
\big( \omega (s_1)\,\omega (s_2)\big)^{1/2} \leq A \,  \omega
\big(\tfrac{s_1+s_2}{2}\big).
\end{equation}
We will obtain convolution estimates in only the dimensions $d=2,3$ and $4$:
\begin{theorem}\label{d=2}
Suppose $d=2$ and assume \eqref{ineq2}. Then there is the Lorentz space estimate
\begin{equation*}%\label{ineq3}
\|\lambda\ast f\|_{L^{3}}\leq C(A)\,\|f\|_{L^{3/2,3}}.
\end{equation*}
\end{theorem}
\begin{theorem}\label{d=3}
Suppose $d=3$ and assume \eqref{ineq2}. Then, for any $\epsilon >0$,  there is the Lorentz space estimate
\begin{equation*}%\label{ineq4.15}
\|\lambda\ast f\|_{L^{2}}\leq C(A)\,\|f\|_{L^{3/2,2-\epsilon}}.
\end{equation*}
\end{theorem}
\begin{theorem}\label{d=4}
Suppose $d=4$ and assume \eqref{ineq2}. If 
\begin{equation}\label{ineq5.15}
\frac{1}{p}-\frac{1}{q}=\frac{1}{10}\ \text{and}\ \frac{4}{10}<\frac{1}{p}<\frac{7}{10}
\end{equation}
then there is the Lebesgue space estimate
\begin{equation*}%\label{ineq3}
\|\lambda\ast f\|_{L^{q}}\leq C(A,p)\,\|f\|_{L^{p}}.
\end{equation*}
\end{theorem}
\noindent Here are some comments:

\noindent(a) Theorem \ref{d=2} is the best possible Lorentz space estimate, 
even in the nondegenerate case $\phi(t)=t^2 /2$. It implies the sharp $L^p \rightarrow L^q$
mapping property, an $L^{3/2}\rightarrow L^3$ estimate.

\noindent(b) Theorem \ref{d=3} is analogous to a result 
from \cite{DLW} for polynomial curves (whose proof we will follow). Theorem \ref{d=3} implies the  
sharp $L^p \rightarrow L^q$ estimates, which hold for 
\begin{equation}\label{ineq4.2}
\frac{1}{p}-\frac{1}{q}=\frac{1}{6}\ \text{and}\ \frac{1}{2}\leq\frac{1}{p}\leq\frac{2}{3}.
\end{equation}
But there are sharper Lorentz space estimates 
for the nondegenerate case $\phi (t)=t^3 /6$ in \cite{BS}
and for polynomial curves (for all 
dimensions $d$) in \cite{Stov2}. 

\noindent(c) Theorem \ref{d=4} is much less satisfactory. One would like, 
for example, at least the sharp $L^p \rightarrow L^q$ mapping properties, which correspond to the endpoints 
in \eqref{ineq5.15}. 

\noindent(d) An analog of Theorem \ref{d=4} for all dimensions $d$, as well as analogs of the 
endpoint results of \cite{Stov1} and \cite{Stov2}, might follow from an analog of the band structure
construction of \cite{C} for the curves and measures considered in this note. But, in view of the
 complicated nature of a Jacobian determinant associated with our operators, it is not clear how to obtain such a band structure.

%\noindent(b) For the sake of comparison we mention that Stovall \cite{Stov2} has proved almost 
%sharp Lorentz space estimates for convolution with affine arclength measure on polynomial curves, with 
%constants depending only on the polynomial's degree. (If $d=2$, the sharp version of Stovall's result 
%follows from \cite{O2} and \cite{BOS1}.) 

\noindent(e) The papers \cite{D}, \cite{Choi1}, and \cite{Choi2} contain
some earlier results for convolution with affine arclength measures in dimensions $2$ and $3$.

Section \ref{proofs} contains the proofs of Theorems 1.2--1.4 and  \S 3 contains proofs
for the lemmas required in \S 2.

\section{Proofs of Theorems}\label{proofs}

\noindent{\bf{Proof of Theorem \ref{d=2}.}}
\begin{proof}
According to the proof of Theorem 5 in \cite{O2}, which abstracts an argument from \cite{BOS1}, it 
is enough to establish the estimate 
\begin{equation*}%\label{ineq4}
\int_a^b \Big(\int_a^b \chi_E \big(\gamma (t_2)-\gamma (t_1)\big)\, \omega(t_2)\, dt_2\Big)^2 \omega (t_1)\, dt_1
\leq C(A)\, |E|
\end{equation*}
for measurable $E\subset \bbR^2$. The inequality 
\begin{equation*}
\int_a^b \Big(\int_{t_1}^b \chi_E \big(\gamma (t_2)-\gamma (t_1)\big)\, \omega(t_2)\, dt_2\Big)^2 
\omega (t_1)\, dt_1 \leq 4\, |E|
\end{equation*}
from \cite{O1}, which is true without any additional hypothesis like \eqref{ineq2}, shows that 
it suffices to establish the estimate
\begin{equation}\label{ineq5}
\int_a^b \Big(\int_a^{t_1} \chi_E \big(\gamma (t_2)-\gamma (t_1)\big)\, 
\omega(t_2)\, dt_2\Big)^2 \omega (t_1)\, dt_1\leq C(A)\, |E|.
\end{equation}
The mapping 
$$
(t_1 ,t_2 )\mapsto \gamma (t_2 )-\gamma (t_1 )
$$ 
is one-to-one by the convexity of the curve $\gamma$. 
If $J(t_1 ,t_2 )$ is the absolute value of the Jacobian determinant of this mapping, then \eqref{ineq5} 
is equivalent to 
\begin{equation}\label{ineq6}
\int \Big(\int \chi_\Omega  (t_1 ,t_2 )\, 
\omega(t_2)\, dt_2\Big)^2 \omega (t_1)\, dt_1\leq C(A)\, \int\int \chi_\Omega  (t_1 ,t_2 )\,
J(t_1 ,t_2 )\, dt_2 \, dt_1
\end{equation}
if $\Omega\subset\{ (t_1 ,t_2 ):a<t_2 <t_1 <b\}$.
We will need the following estimate, a consequence of Lemma \ref{lemma0} below,
\begin{equation}\label{ineq7.5}
J(t_1 ,t_2 )\geq
c(A)\,|t_1 -t_2 |\,\omega(t_2 )\,\omega (t_1)^2. 
\end{equation}

\begin{lemma}\label{lemma0}
Suppose $\gamma$ is as in \eqref{ineq0} and let $J(t_1 ,\dots ,t_d )$ be the absolute value of the 
Jacobian determinant of one of the mappings
$$
(t_1 ,\dots ,t_d )\mapsto \gamma (t_1 )\pm \gamma (t_2 )\pm\cdots\pm\gamma (t_d).
$$
Suppose that \eqref{ineq2} holds and that $n_1 \leq \cdots\leq n_d$ are positive numbers satisfying 
$n_1 +\cdots +n_d =d(d+1)/2$. Suppose that 
$\{i_1 ,\dots ,i_d\}=\{1,\dots ,d\}$ and that 
$a<t_{i_1} <\cdots <t_{i_d} <b$. Then 
\begin{equation}\label{ineq3.0}
J(t_1 ,\dots ,t_d )\ge c\, \Big(\prod_{j=1}^d \omega (t_{i_j} )^{n_j}\Big)  V(t_1 ,\dots ,t_d ),
\end{equation}
where 
$$
V(t_1 ,\dots ,t_d )=\big|
\prod_{1\le i<j\le d}(t_j -t_i )
\big|
$$
and where $c$ depends only on $A$ from \eqref{ineq2} and on $n_1 ,\dots , n_d$.
\end{lemma}

\noindent Given \eqref{ineq7.5}, inequality \eqref{ineq6} will follow from 
\begin{equation}\label{ineq8}
\int \Big(\int \chi_\Omega  (t_1 ,t_2 )\, 
\omega(t_2)\, dt_2\Big)^2 \omega (t_1)\, dt_1\leq 
\end{equation}
\begin{equation*}
C\, \int\Big(\int \chi_\Omega  (t_1 ,t_2 )\,
\omega (t_1 )\,|t_1 -t_2 |\,\omega (t_2)
\, dt_2 \Big)\omega (t_1 )\, dt_1 .
\end{equation*}
To see \eqref{ineq8}, we will use the following lemma.

\begin{lemma}\label{lemma1}
Suppose $\omega$ is nonnegative and nondecreasing on some interval $(c,d]$. Suppose that 
$c_1 ,\dots ,c_K \in\bbR$. For $\rho>0$ let
$$
E_\rho =\{t\in(c,d]: \omega (t)^{K-1} \,\omega (d)\prod_{\ell =1}^K |t-c_l |\leq \rho ^K\}.
$$
Then 
$$
\int_{E_\rho}\omega (t)\,dt\leq C(K)\,\rho .
$$
\end{lemma}
\noindent Indeed, fix $t_1$ and define $\rho$ by
$$
\rho=\frac{1}{2\,C(1)}\int \chi_\Omega  (t_1 ,t_2 )\, \omega(t_2)\, dt_2 ,
$$
where $C(1)$ is the constant in Lemma \ref{lemma1} corresponding to $K=1$. It follows from 
Lemma \ref{lemma1} (with $d=t_1$) that
$$
\int \chi_\Omega  (t_1 ,t_2 )\,
\omega (t_1 )\,|t_1 -t_2 |\,\omega (t_2)\, dt_2 \geq \frac{1}{4\,C(1)} \Big(\int \chi_\Omega  (t_1 ,t_2 )\, 
\omega(t_2)\, dt_2\Big)^2 .
$$
Now integrating with respect to $t_1$ gives \eqref{ineq8}.

\end{proof}

\noindent{\bf{Proof of Theorem \ref{d=4}.}}
\begin{proof}

We will apply the iterated $TT^*$ method introduced by Christ in \cite{C} and 
(see, e.g., the discussion and references in \cite{Stov2})
employed by many others since then.
Thus, assuming some familiarity with Christ's method, Theorem \ref{d=4} will follow if we 
establish the inequality \eqref{ineq9} below, where 
$E$, $\alpha$, and $\beta$ are as follows:
let $\Omega\subset (a,b)^4$ be a set of the form
\begin{equation*}%\label{ineq8.25}
\Omega =\{(t_1 ,t_2 ,t_3 ,t_4):t_1 \in\Omega_0 ,\,t_2 \in \Omega (t_1 ),\, t_3 \in \Omega (t_1 ,t_2 ),\,
\, t_4 \in \Omega (t_1 ,t_2 ,t_3 )\}
\end{equation*}
where 
\begin{equation}\label{ineq8.5}
\lambda (\Omega_0 )\geq \beta >0,\ \lambda \big(\Omega (t_1 ) \big)\geq \alpha >0 \ \text{for each}\  t_1\in\Omega_0 ,
\ 
\end{equation}
\begin{equation}\label{ineq8.6}
\lambda \big(\Omega (t_1 ,t_2) \big)\geq \beta \  \text{whenever}\  
t_1\in\Omega_0 ,\,t_2\in\Omega (t_1),\ \text{and}\ 
\end{equation}
\begin{equation}\label{ineq8.7}
\lambda \big(\Omega (t_1 ,t_2 ,t_3) \big)\geq \alpha \  \text{whenever}\  t_1\in\Omega_0 ,\,t_2\in\Omega (t_1),\ 
 t_3 \in \Omega (t_1 ,t_2 ).
\end{equation}
(Here we are
writing $\lambda$ for the measure 
$$
d\lambda (t)=\omega (t)\, dt=\phi ^{(4)} (t)^{1/10}dt
$$
on $(a,b)$ as well as for its image on $\gamma$.) The set $E$ is defined by
$$
E=\{\gamma (t_1 )-\gamma (t_2 )+\gamma (t_3 )-\gamma (t_4 ):(t_1 ,t_2 ,t_3 ,t_4 )\in \Omega \},
$$
and the desired inequality is
\begin{equation}\label{ineq9}
|E|\geq c(A)\,\alpha^7 \beta^3 .
\end{equation}

%
%the claim above will follow from the inequality
%
%\begin{equation}\label{ineq10}
%\int_{\Omega_0}
%\Big( \int_{\Omega (t_1 )}\Big[
%\int_{\Omega (t_1 ,t_2 )}
%\Big(\int_{\Omega (t_1 ,t_2 ,t_3)}\omega (t_4 )\,dt_4
%\Big)^2
%\omega (t_3 )\, dt_3
%\Big]^{2/3}\omega (t_2 )\, dt_2\Big)^3                            \omega (t_1 )\,dt_1
%\leq 
%\end{equation}
%\begin{equation*}
%C(A)\int_{\Omega}J(t_1 ,t_2 ,t_3 ,t_4)\, dt_4 \, dt_3\,dt_2\,dt_1 .
%\end{equation*}
By passing to a subset of $\Omega$ and replacing $\alpha$ and $\beta$ by $\alpha /24$ and $\beta /24$,
we can assume that there is some permutation $\{i_1 ,i_2 ,i_3 ,i_4\}$ of 
$\{1,2,3,4\}$ such that if $(t_1 ,t_2 ,t_3 ,t_4)\in\Omega$ then 
\begin{equation*}%\label{ineq11}
t_{i_1}<t_{i_2}<t_{i_3}<t_{i_4}. 
\end{equation*}
If $J(t_1 ,t_2 ,t_3 ,t_4)$ is the absolute value of the Jacobian determinant of the mapping
$$
(t_1 ,t_2 ,t_3 ,t_4)\mapsto \gamma (t_1 )-\gamma (t_2 )+\gamma (t_3 )-\gamma (t_4 ),
$$
we will use the following inequality, a consequence of Lemma \ref{lemma0},
\begin{equation}\label{ineq12}
J(t_1 ,t_2 ,t_3 ,t_4 )\geq c(A)\,
\omega (t_{i_4})^9
\Big(\prod_{j=1}^3 
\omega (t_{i_j})\Big)^{1/3}\,V(t_1 ,t_2 ,t_3 ,t_4 ). 
\end{equation}
We will also need the following lemma.
\begin{lemma}\label{lemma2}
Suppose $\omega$ is nonnegative and nondecreasing on an interval $[c,d)$. Suppose $\eta >0$ and $r>1$
satisfy
$$
\eta<\frac{1}{r'}\ \dot=\ 1-\frac{1}{r}.
$$
Suppose $E\subset [c,d)$ and let 
$$
\rho =\int_E \omega (t)\,dt .
$$
Then, for $t_0 \in\bbR$,
\begin{equation*}%\label{ineq13}
\rho^{1+r\eta}\ \omega (c)^{r-(1+r\eta )}\leq C(\eta ,r) \, \int_E \omega (t)^r \, |t-t_0 |^{r\eta}\,dt.
\end{equation*}
If $t_1 ,t_2 \in\bbR$ then also 
\begin{equation*}%\label{ineq14}
\rho^{1+r\eta}\,|t_2 -t_1 |^{r\eta }\,\omega (c)^{r-(1+r\eta )}\leq 
C(\eta ,r) \, \int_E \omega (t)^r \, (|t-t_1 |\cdot|t-t_2 |)^{r\eta}\,dt.
\end{equation*}
\end{lemma}
Now define $I$ by
\begin{equation*}%\label{ineq15}
I=
\int_{\Omega_0} \int_{\Omega (t_1 )}\int_{\Omega (t_1 ,t_2 )}\int_{\Omega (t_1 ,t_2 ,t_3)}
\omega (t_{i_4})^9\,
\Big(\prod_{j=1}^3 
\omega (t_{i_j})\Big)^{1/3}\,V(t_1 ,t_2 ,t_3 ,t_4 )
\,dt_4 \, dt_3 \, dt_2 \,dt_1 
\end{equation*}
so that, because of \eqref{ineq12}, we have 
\begin{equation}\label{ineq15.5}
|E|\geq c(A)\,I.
\end{equation}
(The change of variables needed for the estimate \eqref{ineq15.5} is justified as in \cite{DM2}, p. 549.)
In view of \eqref{ineq15.5}, \eqref{ineq9} will follow if we show that
\begin{equation}\label{15.7}
I\gtrsim \alpha^7 \beta^3.
\end{equation}
(The constants implied by $\lesssim$ and $\gtrsim$ will not depend on any parameters.)

We will, unfortunately, need to consider several cases. To begin, if $4=i_4$, we will use 
Lemma \ref{lemma2} with $r=5$ and $\eta =3/5$ to estimate
\begin{equation*}%\label{ineq16}
\int_{\Omega (t_1 ,t_2 ,t_3)}\omega (t_4 )^5 \,\prod_{j=1}^3 |t_4 -t_j |\,dt_4 \geq 
\int_{\Omega (t_1 ,t_2 ,t_3)}\omega (t_4 )^5 \, |t_4 -t_{i_3} |^3 \,dt_4 \gtrsim
\end{equation*}
\begin{equation*}
\,\Big(\int_{\Omega (t_1 ,t_2 ,t_3)}\omega (t_4 )\,  dt_4 \Big)^4\omega (t_{i_3}).
\end{equation*}
With the inequality $\omega (t_{i_4})\ge\omega (t_{i_3})$ this gives 
\begin{equation}\label{ineq17}
I\gtrsim
\int_{\Omega_0} \int_{\Omega (t_1 )}\int_{\Omega (t_1 ,t_2 )}
\Big(\int_{\Omega (t_1 ,t_2 ,t_3)}\omega (t_4 )\,dt_4 \Big)^4 \cdot  
\end{equation}
\begin{equation*}
\omega (t_{i_3})^{16/3}\,\omega (t_{i_2})^{1/3}\,\omega (t_{i_1})^{1/3}\,V(t_1 ,t_2 ,t_3 )\,dt_3 \,dt_2 \,dt_1 .
\end{equation*}
If $4 =i_{k_0}$ for some $k_0 =1,2,3$, then 
$$
\omega (t_4 )^{1/3}\,\omega (t_{i_4})^{11/3}\ge \omega (t_4)^3 \,\omega (t_{i_{k_0 +1}})
$$
by the monotonicity of $\omega$. Thus 
\begin{equation*}
\int_{\Omega (t_1 ,t_2 ,t_3)}\omega (t_4 )^{1/3}\,\omega (t_{i_4})^{11/3}\prod_{j=1}^3 |t_4 -t_j |\,dt_4 \geq 
\int_{\Omega (t_1 ,t_2 ,t_3)}\omega (t_4)^3\,\omega (t_{i_{k_0 +1}})
\prod_{j=1}^3 |t_4 -t_j |\,dt_4 \gtrsim 
\end{equation*}
\begin{equation*}
\Big(\int_{\Omega (t_1 ,t_2 ,t_3)}\omega (t_4 )\,dt_4 \Big)^4 ,
\end{equation*}
where the last inequality follows from an application of Lemma \ref{lemma1} as at the end of 
the proof of Theorem \ref{d=2} but with $K=3$ instead of $K=1$. Therefore 
\begin{equation}\label{ineq18}
I\gtrsim
\int_{\Omega_0} \int_{\Omega (t_1 )}\int_{\Omega (t_1 ,t_2 )}
\Big(\int_{\Omega (t_1 ,t_2 ,t_3)}\omega (t_4 )\,  dt_4 \Big)^4 \times
\end{equation}
\begin{equation*}
\omega (t_{i_4})^{16/3}\,
\prod_{k=1,k\not=k_0}^3 \omega (t_{i_k})^{1/3}\,
V(t_1 ,t_2 ,t_3 )\,dt_3 \,dt_2 \,dt_1 .
\end{equation*}
Now if $\{j_1 ,j_2 ,j_3 \}$ is the permutation of $\{1,2,3\}$ such that $t_{j_1}<t_{j_2}<t_{j_3}$ whenever 
$t_1 \in\Omega_0 ,\,t_2 \in \Omega (t_1),\, t_3 \in\Omega (t_1 ,t_2 )$, then \eqref{ineq17},  \eqref{ineq18},
and \eqref{ineq8.7} imply that
\begin{equation}\label{ineq19}
I\gtrsim 
\alpha^4 \int_{\Omega_0} \int_{\Omega (t_1 )}\int_{\Omega (t_1 ,t_2 )}
\omega (t_{j_3})^{16/3}\,\omega (t_{j_2})^{1/3}\,\omega (t_{j_1})^{1/3}\,V(t_1 ,t_2 ,t_3 )\,dt_3 \,dt_2 \,dt_1 .
\end{equation}

If $3=j_3$,
\begin{equation*}%\label{ineq20}
\int_{\Omega (t_1 ,t_2 )} \omega (t_3 )^3 \,|t_3 -t_2 |\cdot |t_3 -t_1 |\,dt_3 \geq
|t_1 -t_2 |\int_{\Omega (t_1 ,t_2 )} \omega (t_3 )^3 \,|t_3 -t_{j_2} |\,dt_3 \gtrsim
\end{equation*}
\begin{equation*}
|t_1 -t_2 |\Big(\int_{\Omega (t_1 ,t_2 )} \omega (t_3 )\, dt_3 \Big)^2 \omega (t_{j_2} ),
\end{equation*}
where the $\gtrsim$ results from an application of Lemma \ref{lemma2} with $r=3$ and $\eta =1/3$. 
With \eqref{ineq19}, \eqref{ineq8.6}, and the monotonicity of $\omega$ this gives 
\begin{equation}\label{ineq21}
I\gtrsim 
\alpha^4 \beta^2 \int_{\Omega_0} \int_{\Omega (t_1 )}
\omega (t_{j_2})^{7/2}\,\omega (t_{j_1})^{1/2}\,|t_1 -t_2 |^2 \,dt_2 \,dt_1 .
\end{equation}
If $3=j_2$, the second conclusion of Lemma \ref{lemma2}, with $r=13/6$ and $\eta =6/13$, gives 
$$
\int_{\Omega (t_1 ,t_2 )} \omega (t_3 )^{13/6}|t_3 -t_{2} |\cdot |t_3 -t_{1} |\,dt_3 \gtrsim
|t_{1}-t_{2}|\Big(\int_{\Omega (t_1 ,t_2 )} \omega (t_3 )\, dt_3 \Big)^2 \omega (t_{j_1})^{1/6}.
$$
From \eqref{ineq19} it then follows that 
\begin{equation}\label{ineq22}
I\gtrsim 
\alpha^4 \beta^2 \int_{\Omega_0} \int_{\Omega (t_1 )}
\omega (t_{j_3})^{7/2}\,\omega (t_{j_1})^{1/2}\,|t_1 -t_2 |^2 \,dt_2 \,dt_1 .
\end{equation}
And if $3=j_1$ then 
$$
\int_{\Omega (t_1 ,t_2 )} \omega (t_3 )\omega (t_{j_2})\,|t_3 -t_{2} |\cdot |t_3 -t_{1} |\,dt_3 \geq
|t_{1}-t_{2}|\int_{\Omega (t_1 ,t_2 )} \omega (t_3 )\, \omega (t_{j_2})\,|t_3 -t_{j_2} |dt_3 \gtrsim
$$
$$
|t_{1}-t_{2}|\Big(\int_{\Omega (t_1 ,t_2 )} \omega (t_3 )\, dt_3 \Big)^2
$$
by Lemma \ref{lemma1} with $K=1$, and so  \eqref{ineq19} gives 
\begin{equation}\label{ineq23}
I\gtrsim 
\alpha^4 \beta^2 \int_{\Omega_0} \int_{\Omega (t_1 )}
\omega (t_{j_3})^{7/2}\,\omega (t_{j_2})^{1/2}\,|t_1 -t_2 |^2 \,dt_2 \,dt_1 .
\end{equation}

Thus if $\{k_1 ,k_2  \}$ is the permutation of $\{1,2\}$ such that $t_{k_1}<t_{k_2}$ whenever 
$t_1 \in\Omega_0 ,\,t_2 \in \Omega (t_1)$, then \eqref{15.7} will follow from 
\eqref{ineq21}, \eqref{ineq22}, and \eqref{ineq23} if we establish that
\begin{equation}\label{ineq24}
\int_{\Omega_0} \int_{\Omega (t_1 )}\omega (t_{k_2})^{7/2}\,\omega (t_{k_1})^{1/2}\,|t_1 -t_2 |^2 \,dt_2 \, dt_1\gtrsim
\alpha^3 \beta.
\end{equation}
If $2=k_2$, then 
$$
\int_{\Omega (t_1 )}\omega (t_{2})^{7/2}\,|t_1 -t_2 |^2 \,dt_2 \gtrsim 
\Big(
\int_{\Omega (t_1 )}\omega (t_{2}) \,dt_2
\Big)^3 \omega (t_1 )^{1/2}
$$
by Lemma \ref{lemma2} with $r=7/2$, $\eta =4/7$, and \eqref{ineq24} follows from \eqref{ineq8.5}.
If $2=k_1$, then 
$$
\omega (t_{k_2})^{7/2}\,\omega (t_{k_1})^{1/2}\ge \omega (t_1 )^2 \omega (t_2 )^2
$$
and
$$
\int_{\Omega (t_1 )}\omega (t_{2})^{2}\,\omega (t_1 )\,|t_1 -t_2 |^2 \,dt_2 \gtrsim 
\Big(
\int_{\Omega (t_1 )}\omega (t_{2}) \,dt_2
\Big)^3
$$
by Lemma \ref{lemma1} with $K=2$. Again, \eqref{ineq24} follows from \eqref{ineq8.5},
and the proof of Theorem \ref{d=4} is complete.

\end{proof}

\noindent{\bf{Proof of Theorem \ref{d=3}.}}
\begin{proof}
The sharp $L^p \rightarrow L^q$ estimates (for the indices in \eqref{ineq4.2}) can be obtained by the method
of \cite{O}. But to obtain the Lorentz space estimates in Theorem \ref{d=3},   
we will follow the proof of the $d=3$ case in \cite{DLW}, again using the method of Christ. Thus we will begin by establishing the following claim (which, by itself, implies the almost sharp Lebesgue space estimates corresponding to strict inequality in \eqref{ineq4.2}): suppose that 
$\Omega\subset (a,b)^3$ is a set of the form
\begin{equation*}%\label{ineq3.1}
\Omega =\{(t_1 ,t_2 ,t_3 ):t_1 \in\Omega_0 ,\,t_2 \in \Omega (t_1 ),\, t_3 \in \Omega (t_1 ,t_2 )\}
\end{equation*}
where 
\begin{equation}\label{ineq3.002}
\lambda (\Omega_0 )\geq \alpha >0,\ \lambda (\Omega (t_1 ) )\geq \beta >0 \ \text{for each}\  t_1\in\Omega_0 
\ \text{and} \ 
\end{equation}
\begin{equation*}
\lambda (\Omega (t_1 ,t_2) )\geq \alpha \  \text{whenever}\  t_1\in\Omega_0 ,\,t_2\in\Omega (t_1).
\end{equation*}
If 
$$
E=\{\gamma (t_1 )-\gamma (t_2 )+\gamma (t_3 ):(t_1 ,t_2 ,t_3 )\in \Omega \},
$$
then we have 
\begin{equation}\label{ineq3.003}
|E|\geq c(A)\,\alpha^4 \beta^2 .
\end{equation}
As before, 
we can assume that there is some permutation $\{i_1 ,i_2 ,i_3 \}$ of 
$\{1,2,3\}$ such that if $(t_1 ,t_2 ,t_3 )\in\Omega$ then 
\begin{equation*}
t_{i_1}<t_{i_2}<t_{i_3}. 
\end{equation*}
With $J(t_1 ,t_2 ,t_3 )$ the absolute value of the Jacobian determinant of the mapping
$$
(t_1 ,t_2 ,t_3 )\mapsto \gamma (t_1 )-\gamma (t_2 )+\gamma (t_3 ),
$$
we will need the following consequence of Lemma \ref{lemma0}: 
\begin{equation}\label{ineq3.005}
J(t_1 ,t_2 ,t_3  )\geq c(A)\,
\omega (t_{i_3})^5
\Big(\prod_{j=1}^2 
\omega (t_{i_j})\Big)^{1/2}\,V(t_1 ,t_2 ,t_3 ). 
\end{equation}
Define $I$ by
\begin{equation*}%\label{ineq3.6}
I=
\int_{\Omega_0} \int_{\Omega (t_1 )}\int_{\Omega (t_1 ,t_2 )}
\omega (t_{i_3})^5\,
\Big(\prod_{j=1}^2 
\omega (t_{i_j})\Big)^{1/2}\,V(t_1 ,t_2 ,t_3 )
 \, dt_3 \, dt_2 \,dt_1 
\end{equation*}
so that, because of \eqref{ineq3.005}, we have 
\begin{equation*}%\label{ineq3.7}
|E|\geq c(A)\,I.
\end{equation*}
(Again, the change of variables here is justified as in \cite{DM2}.)
Then \eqref{ineq3.003} will follow from
\begin{equation}\label{ineq3.008}
I\gtrsim \alpha^4\, \beta^2.
\end{equation}
Since the proof of \eqref{ineq3.008} is very similar to the proof of \eqref{15.7}, we will only sketch the argument. 
The first step is to obtain the inequality
\begin{equation*}%\label{ineq3.9}
I\gtrsim \alpha^3 \int_{\Omega_0}\int_{\Omega (t_1 )}
\omega (t_{j_2})^{5/2}\omega (t_{j_1})^{1/2}\,|t_1 -t_2 |\,dt_2 \, dt_1,
\end{equation*}
where $\{t_1 ,t_2 \}=\{ t_{j_1},t_{j_2}\}$ and $ t_{j_1}<t_{j_2}$. Recalling \eqref{ineq3.002}, this is done by using  
 Lemma \ref{lemma2},  
with $r=7/2$ and $\eta =4/7$, if $3=i_3$ and by using Lemma \ref{lemma1}, with $K=2$,  if $3=i_2$ or $3=i_1$.
The proof of \eqref{ineq3.008} is then concluded by showing that 
\begin{equation*}
 \int_{\Omega_0}\int_{\Omega (t_1 )}
\omega (t_{j_2})^{5/2}\omega (t_{j_1})^{1/2}\,|t_1 -t_2 |\,dt_2 \, dt_1
\gtrsim \beta^2 \,\alpha 
\end{equation*}
by using Lemma \ref{lemma2} with $r=5/2$, $\eta =2/5$ if $t_1 <t_2$ and Lemma \ref{lemma1} with $K=1$ if $t_2 <t_1$.
This proves \eqref{ineq3.008} and thus, as mentioned above, establishes the almost-sharp
Lebesgue space bounds by the method of 
\cite{C}. 

To obtain the Lorentz space bounds claimed in Theorem \ref{d=3}, we follow the proof 
of the analogous result in \cite{DLW} (itself based on a further argument of Christ \cite{C2}). Thus it is enough
to establish an analogue of Lemma 1 in \cite{DLW} for our curves $\gamma$. The crux of the matter is to show
the following: if 
$\Omega\subset (a,b)^3$ is a set of the form
\begin{equation*}
\Omega =\{(t_1 ,t_2 ,t_3 ):t_1 \in\Omega_0 ,\,t_2 \in \Omega (t_1 ),\, t_3 \in \Omega (t_1 ,t_2 )\},
\end{equation*}
where 
\begin{equation*}
\lambda (\Omega_0 )\geq \beta >0,\ \lambda (\Omega (t_1 ) )\geq \beta |G|/|E| >0 \ \text{for each}\  t_1\in\Omega_0 
\ \text{and} \ 
\end{equation*}
\begin{equation*}
\lambda (\Omega (t_1 ,t_2) )\geq \delta >0 \  \text{whenever}\  t_1\in\Omega_0 ,\,t_2\in\Omega (t_1),
\end{equation*}
and if 
$$
E'=\{\gamma (t_1 )-\gamma (t_2 )+\gamma (t_3 ):(t_1 ,t_2 ,t_3 )\in \Omega \},
$$
then we have 
\begin{equation*}
|E'|\geq c(A)\,\delta^3 \Big(\frac{\beta |G|}{|E|}\Big)^2\, \beta .
\end{equation*}
This can be established by exactly the argument  
given above for \eqref{ineq3.003}.

\end{proof}

\section{Proofs of lemmas}\label{lemmas}

\noindent{\bf{Proof of Lemma \ref{lemma0}.}}
\begin{proof}
Assume without loss of generality that $a<t_1 <\dots <t_d <b$.
It is enough to prove the lemma in the special case when each $n_j$ can be written  
\begin{equation}\label{ineq3.05}
n_j =\frac{d(d+1)}{2}\cdot\frac{l_j}{2^n}
\end{equation}
for some large integer $n$ and positive integers $l_j$.
(To see this, find $n$ and $l_1 \leq \cdots \leq l_d$ such that 
$$
\sum_{j=1}^d l_j =2^n \ \text{and}\ 
n_j \le\frac{d(d+1)}{2}\cdot\frac{l_j}{2^n},\,j=2,\dots ,d.
$$
Then note that 
$$
\prod_{j=1}^d \omega (t_j )^{n_j}\leq \prod_{j=1}^d \omega (t_j )^{\tfrac{d(d+1)l_j}{2^{n+1}}}
$$
by the monotonicity of $\omega$.)

It is shown in \cite{BOS2} that there exists a nonnegative function $\psi =\psi (u;t_1 ,\dots ,t_d )$ supported in 
$[t_1 ,t_d ]$ such that 
\begin{equation}\label{ineq3.1}
J(t_1 ,\dots ,t_d )=\int_{t_1}^{t_d}\omega (u)^{d(d+1)/2}\,\psi (u)\,du .
\end{equation}
The choice $\phi (t)=t^d /d!$ in \eqref{ineq0} shows that 
\begin{equation*}%\label{ineq3.2}
\int_{t_1}^{t_d}\psi (u)\,du =c(d)\, V(t_1 ,\dots ,t_d ).
\end{equation*}
For $\delta >0$, define
$$
t_\delta =\big(1-(d-1)\delta\big)t_d +\delta (t_1 +t_2 +\cdots +t_{d-1}).
$$
The inequality \eqref{ineq3.0} will follow from 
\eqref{ineq3.1}, the monotonicity of $\omega$, the inequality
\begin{equation}\label{ineq3.3}
\int_{t_\delta}^{t_d}\psi (u; t_1  ,
\dots ,t_d )\,du \ge c(\delta )\, V(t_1 ,\dots ,t_d ),
\end{equation}
and the fact that there is a $\delta =\delta (n_1 ,\dots n_d )>0$ such that
\begin{equation}\label{ineq3.4}
\omega (t_\delta )^{d(d+1)/2}\geq c(A;n_1 ,\dots ,n_d )\,\prod_{j=1}^d \omega (t_j )^{n_j}.
\end{equation}

The proof of \eqref{ineq3.3} is by induction on $d$. Since 
$$
\psi (u;t_1 ,t_2 )=\chi_{[t_1 ,t_2 ]}(u),
$$ 
the case $d=2$ is clear. The inductive step requires an identity from \cite{BOS2}:
\begin{equation*}%\label{ineq3.5}
\psi (u;t_1 ,\dots ,t_d )=
\int_{t_1}^{t_2}\cdots \int_{t_{d-1}}^{t_d}\psi (u;s_1 ,\dots ,s_{d-1})\,ds_1 \cdots ds_{d-1}. 
\end{equation*}
Thus 
\begin{equation}\label{ineq3.6}
\int_{t_\delta}^{t_d}\psi (u)\,du =
\int_{t_1}^{t_2}\cdots \int_{t_{d-1}}^{t_d}
\int_{\{u\ge t_\delta\}}\psi (u;s_1 ,\dots ,s_{d-1})\,du\,ds_1 \cdots ds_{d-1}.
\end{equation}
We need the following additional fact from \cite{BOS2}: suppose $\lambda_j \in (0,1)$ for 
$j=1 ,\dots ,d-1$ and let 
$$
t'_j= \lambda_{j}t_{j}+(1-\lambda_j )t_{j+1}
$$
for $j=1 ,\dots ,d-1$.
Then
\begin{equation}\label{ineq3.7}
\int_{t'_1}^{t_2}\cdots \int_{t'_{d-1}}^{t_d}V(s_1 ,\dots ,s_{d-1})\,ds_1 \cdots ds_{d-1}
\ge c(\lambda_1 ,\dots ,\lambda_{d-1})\, V(t_1 ,\dots ,t_d ).
\end{equation}
Now choose $\lambda_1 ,\dots ,\lambda_{d-1}\in(0,1)$ and $\delta' >0$ such that
if $ t'_j \leq s_j \leq t_{j+1}$ for $j=1,\dots,d-1$ then 
\begin{equation}\label{ineq3.75}
s_{\delta'} \doteq \big(1-(d-2)\delta' \big)s_{d-1} +\delta' (s_1 +s_2 +\cdots +s_{d-2})\ge
\end{equation}
$$
t_\delta =\big(1-(d-1)\delta\big)t_d +\delta (t_1 +t_2 +\cdots +t_{d-1}).
$$
(Here is how to make this choice: we can assume that $t_d =1$.
%
%\begin{equation}\label{ineq3.8}
%0<t_1 <\cdots <t_d =1.
%\end{equation}
If $$1-(d-2)\delta ' >0,$$ then \eqref{ineq3.75} holds for all 
$s_j \in [ t'_j ,t_{j+1}]$ if and only if it holds for $s_j =t'_j$. So fix $s_j =t'_j$ for $j=1,\dots ,d-1$.
Then, with $\lambda =(\lambda _1 ,\dots ,\lambda_{d-1})$,
$$
s_{\delta '}=(1-\lambda_{d-1})\big( 1-(d-2)\delta' \big)+\sum_{j=1}^{d-1} c_j (\delta ' ,\lambda )t_j
$$
where 
$$
|c_j (\delta ' ,\lambda )|=O(\delta ' +\|\lambda \|).
$$
Assume that $\delta '$ and $\lambda$ are chosen so that 
$$
(1-\lambda_{d-1})\big( 1-(d-2)\delta' \big)\ge \big( 1-(d-1)\delta \big)\ \text{and}\ |c_j (\delta ' ,\lambda )|\le \delta .
$$
Since
$$
t_\delta =\big( 1-(d-1)\delta )+\delta (t_1 +\cdots +t_{d-1} ),
$$
it then follows from the fact that $s_{\delta '}=t_\delta =1$ when $t_1 =\cdots =t_{d-1}=1$ that $s_{\delta '}\ge t_\delta$
if $0\le t_j \le1$.)

\noindent Now
\begin{equation*}
\eqref{ineq3.6}
\ge 
\int_{t'_1}^{t_2}\cdots \int_{t'_{d-1}}^{t_d}
\int_{\{u\ge t_\delta\}}\psi (u;s_1 ,\dots ,s_{d-1})\,du\,ds_1 \cdots ds_{d-1}\ge
\end{equation*}
\begin{equation*}
\int_{t'_1}^{t_2}\cdots \int_{t'_{d-1}}^{t_d}
\int_{\{u\ge s_{\delta'}\}}\psi (u;s_1 ,\dots ,s_{d-1})\,du\,ds_1 \cdots ds_{d-1}\ge
\end{equation*}
\begin{equation*}
c (\delta' )\int_{t'_1}^{t_2}\cdots \int_{t'_{d-1}}^{t_d}
V(s_1 ,\dots ,s_{d-1})
\,ds_1 \cdots ds_{d-1}\ge 
\end{equation*}
\begin{equation*}
c(\delta' ;\lambda_1 ,\dots ,\lambda_{d-1})\,V(t_1 ,\dots ,t_d ),
\end{equation*}
where the second inequality is due to \eqref{ineq3.75} and the fact that $\psi (u;s_1 ,\dots ,s_d )$ is nonnegative, the third to the induction hypothesis, and the fourth
to \eqref{ineq3.7}. This completes the proof by induction on $d$ of \eqref{ineq3.3}.

To see \eqref{ineq3.4}, recall from \eqref{ineq3.05} that 
\begin{equation*}
%\label{ineq3.05}
n_j =\frac{d(d+1)}{2}\cdot\frac{l_j}{2^n}
\end{equation*}
for some large integer $n$ and positive integers $l_j$
satisfying
$$
\sum_{j=1}^d \frac{l_j}{2^n}=1.
$$
Choose $\delta >0$ so small that 
$$
\delta< \frac{l_j}{2^{n}}
$$
for $j=1,\dots ,d-1$. Note that, since $t_j <t_d $ if $j<d$,  
$$
t_\delta =\big(1-(d-1)\delta \big)t_d +\delta (t_2 +\cdots +t_{d-1}) \ge \sum_{j=1}^d   \frac{l_j}{2^{n}}\,t_j .
$$
Now the inequality
\begin{equation*}
%\label{ineq7}
A^n \, \omega
\big(\tfrac{s_1+\dots+s_{2^n}}{2^n}\big)\ge
\Big(\prod_{j=1}^{2^n} \omega(s_j)\Big)^{1/2^n} 
\end{equation*}
(which follows from iterating \eqref{ineq2}) 
and the monotonicity of $\omega$ imply that 
$$
\omega (t_\delta )\ge \omega \Big(\frac{1}{2^n}\sum_{j=1}^d   {l_j}\,t_j \Big)
\ge A^{-n}\,\prod_{j=1}^d \omega (t_j )^{\frac{l_j}{2^{n}}}.
$$
This give \eqref{ineq3.4}.

\end{proof}

\noindent{\bf{Proof of Lemma \ref{lemma1}.}}
\begin{proof}
By scaling we can assume that $\rho =1$.
Partition $(c,d]$ into disjoint intervals $I_j =(a_j ,a_{j+1}]$ such that $2^{j}\le \omega \le 2^{j+1}$ on $I_j$.
Assume $d\in I_{j_0}$.
We will need the inequality
\begin{equation}\label{ineq3.9}
|\{t\in\bbR :\prod_{l=1}^K |t-c_l |\leq \tau \}|\leq C(K)\, \tau ^{1/K}, \ \tau >0.
\end{equation}
(To see \eqref{ineq3.9}, observe that $\bbR$ can be partitioned into at most $2K$ intervals $J_p$
with the property that
$$
\prod_{l=1}^K |t-c_l |\geq |t-c_{l(p)}|^K ,\ t\in J_p .)
$$
From \eqref{ineq3.9} it follows that if 
$$
E_j =\{t\in I_j :\omega (t)^{K-1}\,\omega (d)\, \prod_{l=1}^K |t-c_l |\leq 1\},
$$
then
\begin{equation*}%\label{ineq3.10}
|E_j | 
\leq \frac{C(K)}
{2^{\big( j(K-1)+j_0 \big)/K}}.
\end{equation*}
Thus
$$
\int_{E_j} \omega (t)\, dt\le C(K)\, 2^{(j-j_0 )/K},
$$
and the conclusion of Lemma \ref{lemma1} follows by summing a geometric series.

\end{proof}

\noindent{\bf{Proof of Lemma \ref{lemma2}.}}
\begin{proof}
We begin by observing that
\begin{equation*}%\label{ineq3.11}
\rho =\int_E \omega (t)\,|t-t_0 |^{\eta }\,|t-t_0 |^{-\eta }\,dt \leq
\end{equation*}
\begin{equation*}
\Big(\int_E  \omega (t)^r |t-t_0 |^{r \eta}\,dt\Big)^{1/r}
\Big(\int_E   |t-t_0 |^{-r' \eta}\,dt\Big)^{1/r'}\le
\end{equation*}
\begin{equation*}
C(r,\eta )\,\Big(\int_E  \omega (t)^r |t-t_0 |^{r \eta}\,dt\Big)^{1/r}
\, |E|^{1-\eta -1/r}.
\end{equation*}
Thus, by the monotonicity of $\omega$, 
\begin{equation*}
\rho^{1+r\eta }\big(\omega (c)\,|E|\big)^{r-1-r\eta}\leq
\rho^{1+r\eta }\,\rho^{r-1-r\eta}=\rho ^r \leq
\end{equation*}
\begin{equation*}
C(r,\eta )\,\int_E  \omega (t)^r |t-t_0 |^{r \eta}\,dt
\cdot |E|^{r-1-r\eta}.
\end{equation*}
This gives the first conclusion of Lemma \ref{lemma2}.
Using the estimate
$$
\int_E \big(|t-t_1 |\cdot |t-t_2 |\big)^{-r' \eta }\,dt\leq C(r,\eta )\,|E|^{1-r' \eta}\,|t_1 -t_2 |^{-r' \eta},
$$
the second conclusion follows similarly.
 
\end{proof}

\end{document}